# Finite element stress analysis of a combined stacker-reclaimer machine: A design audit report

**Erfan Khodabandeh**[1,2*], **Shahaboddin Amini** [1], **Aliakbar Taghipour** [1]

[1] Central Technical Office, Goharzamin Iron Ore Company, Sirjan, Iran
[2] Mechanical Engineering Department, Amirkabir University of Technology, Tehran, Iran

June 2021

* Corresponding author: Erfan Khodabandeh, Email: khodabandeh.e.1989@gmail.com , Tel: 0098-9132083909

**Abstract**

Design audit or design verification is an important step in engineering of heavy mobile materials handling equipment. Usually, the costumers employ third parties for audition of contractors engineering. Here a part of design audit of a combined stacker-reclaimer machine is reported. This equipment is designed and constructed by a local supplier in Iran for the iron ore pelletizing plants at Goharzamin Iron Ore Company. The structure plays an important role in mobile material handling machines such as Stackers and Reclaimers and its failure and damage may cause considerable financial and human life losses. In this report, the stacker-reclaimer machine's undercarriage including gantry and traveling system are numerically analyzed. The Finite Element Method (FEM) is used for stress prediction under the critical operating loads according to the design standards. The critical areas of the undercarriage are identified and it is observed that, the maximum stress is in the safe range.

**Introduction**

Handling of bulk materials plays a significant role in the modern economy and the volume of handled bulk materials has been steadily increased. Hence various types of bulk handling equipment and trans-shipping procedures have been developed in order to make long-distance transportation more economical. The handling of raw bulk materials constitutes a significant cost portion of the final product. Numerous systems of this kind are installed at mineral processing, coal power plants and many other industries where large volumes of bulk materials need to be handled [1-3].

The materials handling systems are used to stack coal, iron ore, fertilizer, and other bulk materials in piles, or reclaim and move them from storage area to other locations. The operation of unloading from the conveyors and building a pile is called "stacking". The "reclaiming" function is to recover bulk material from a stockpile and transferring to a conveyor line. Stackers and reclaimers are specialized for their respective operations, and a combined stacker-reclaimer is capable of doing both operations [1-3], see Fig. 1.

There are many manufacturers and vendors throughout the world that supply mobile materials handling equipment. According to the design standards [4-10], they need to predict, analyze, and calculate the maximum permissible stresses of the structures. For this purpose, the FEM[1] software packages are widely used [4-11].

The importance of FEA[2] lies on its capability to work easily with complex structural designs and offer insight into the structural behavior and robustness of the designs. Using FEA, engineers gain valuable information on weak points, probable failures and improper design features.

Design audition acts as a quality assurance tool. It reviews the project life cycle evaluating the results yielded during the different stages, from the design phase to implementation. The importance of project (design) audit is highly recommended in the standard [9,10] to prevent machine (design) failure [12,13]. In this report, the calculation and design of a combined stacker-reclaimer machine of the iron ore pelletizing plants at Goharzamin Iron Ore Company is reviewed. This equipment is designed and constructed by a local supplier for the first time in Iran. It has a design capacity of 1800 ton/h in stacking and 2000 ton/h in reclaiming. Relying on numerical FEA in Ansys [14], the two parts of the stacker-reclaimer machine (gantry and traveling system), which withstand the highest stresses, are analyzed. The main goal of this study is to investigate the behavior of these parts under critical loading conditions, and also to identify the critical areas. The maximum stresses in different parts of the structures are compared with permissible values.

**Geometry**

The bucket wheel Stacker Reclaimer and the structural parts are shown in Figs. 1 and 2 respectively. The machine specifications are provided in Table 1 and the structural parts are introduced in Table 2. The machine structure is mainly divided into two parts, the undercarriage and the super structure. The undercarriage structure consists of traveling system and gantry, see Fig. 3. The superstructure is composed of two booms (main and counterweight), a bucket wheel, a slewing platform, A-frames and tension links. In this study, FEA of the undercarriage structure is carried out. The gantry and traveling system are mainly made of ST52 with 355 MPa yield strength. According to the FEM 2.131 standard [8, 9], the equipment is classified in P3 and E5 groups; so, the permissible stress is 240 MPa.

---

[1] Finite Element Method (FEM)
[2] Finite Element Analysis (FEA)

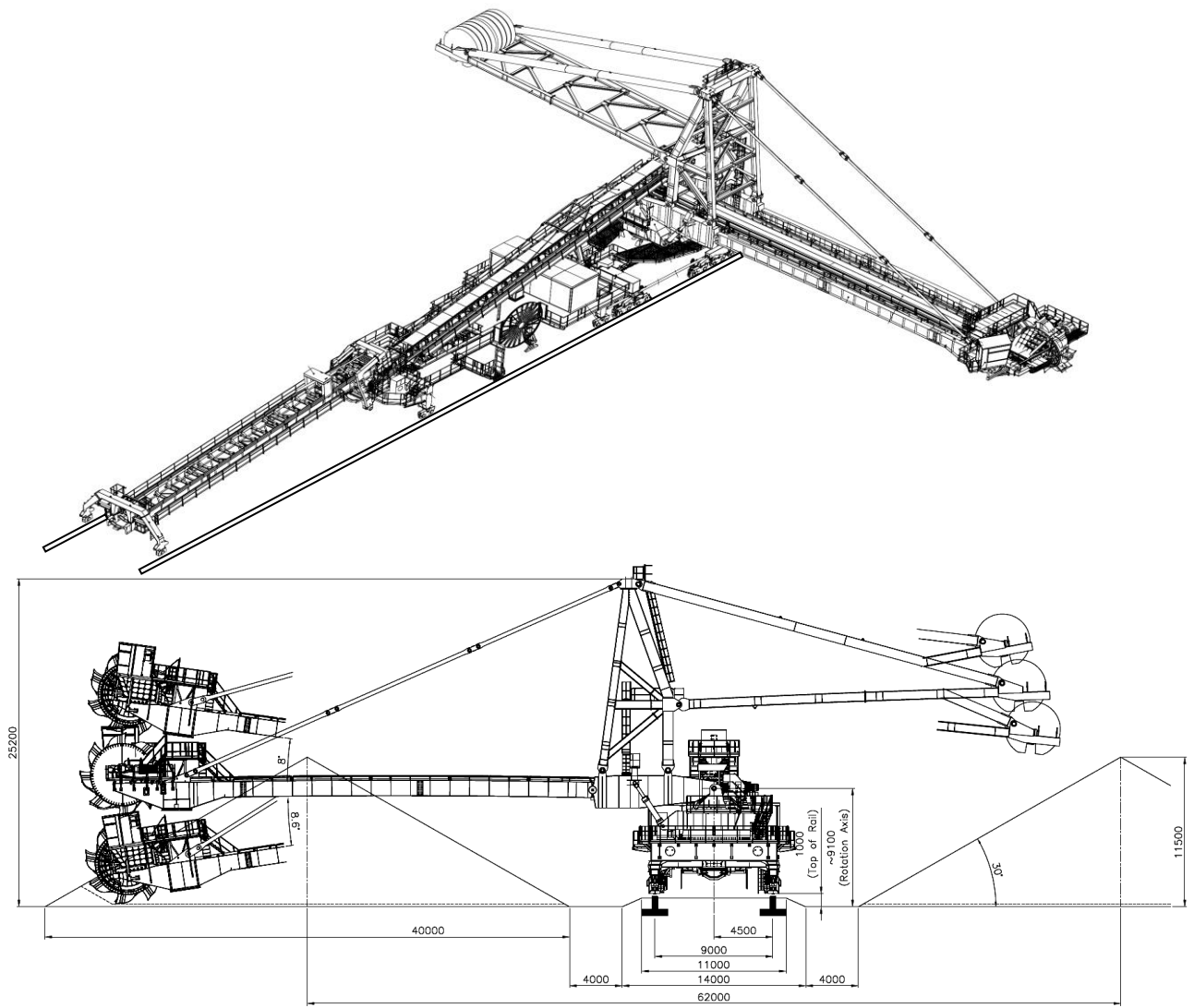

Figure 1. The combined Stacker-Reclaimer at Goharzamin (top: 3D view, bottom: operation in the stock yard).

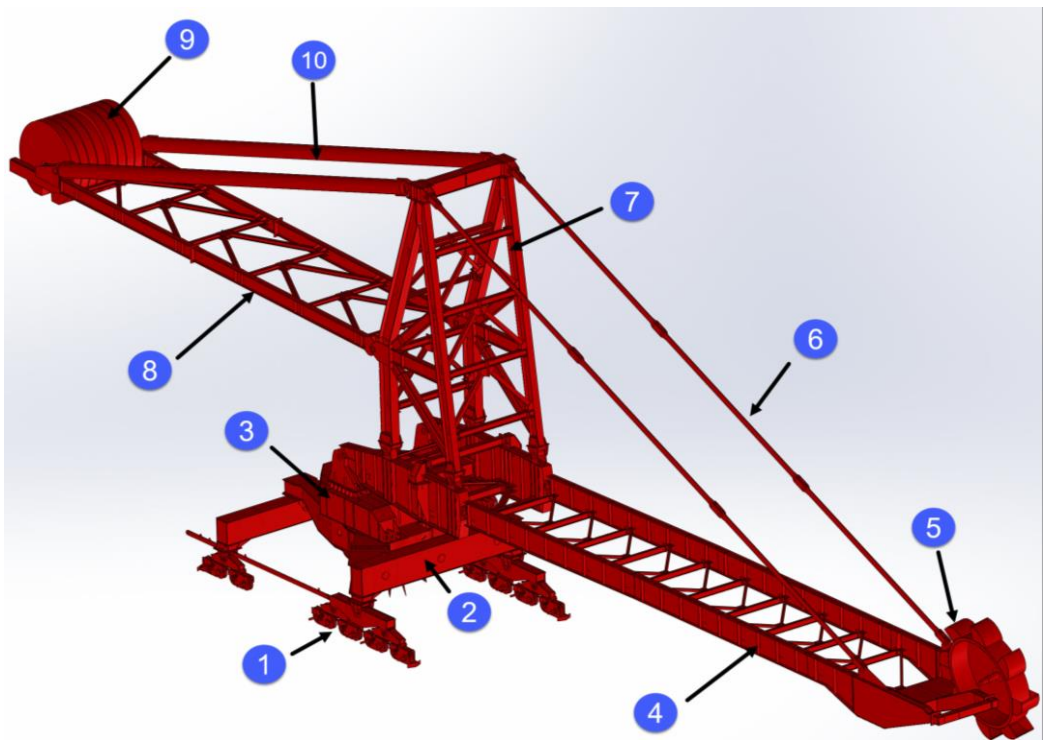

Figure 2. Main structural parts of combined Stacker-Reclaimer (numbered parts are described in Table 2).

Table 1. The combined Stacker-Reclaimer at Goharzamin pelletizing plants.

| Reclaiming (ton/h) | Stacking (ton/h) | Length (m) | Height (m) | Width (m) | Slewing (degree) | Luffing (degree) | Boom (m) | Buckets | Weight (ton) |
|---|---|---|---|---|---|---|---|---|---|
| 2000 | 1800 | 110 | 25 | 12 | -110 to +110 | -9 to +9 | 45 | 8 | 600 |

Table 2. Main structural parts of combined Stacker-Reclaimer.

| No. | Structural parts | | Description | Weight (ton) |
|---|---|---|---|---|
| 1 | Traveling system | Under Carriage | It consists of moving and non-moving wheels, bogeys, primary and secondary balance links, supports for connection to the gantry structure, pins for connecting various parts of the structure together, tension links and other standard parts such as nut, bolts and washers, etc. This machine travel on a rail between stockpiles in the stockyard | 36 |
| 2 | Gantry | | It consists of the portal bogey, and the structure connected to the gantry, supports, gear thrust bearing, connecting pins, etc. | 63 |
| 3 | Slewing platform | Super Structure | It consists of the turning device, supports of the main boom, and other standard parts, including nuts, bolts and washers, etc. | 50 |
| 4 | Main Boom | | It consists of the end, mid and front frames, reinforcement and truss beams, connecting pins, and other standard parts such as nuts, bolts and washers, etc. | 70 |
| 5 | Bucket Wheel | | The bucket wheel unit consists of the rotary wheel body with buckets. | 13.2 |
| 6, 7 | Tension links and A-frame | | They consist of the main frame, tension links of the main boom and the balance boom, connecting pins and other standard parts such as nuts, bolts and washers, etc. | 71 |
| 8, 9, 10 | Counter weight, beam and tension links | | This structure includes longitudinal beams, reinforcement, and truss transverse beams, connecting pins and other standard parts such as nuts, bolts and washers, etc. | 30 |
| | | | Total | 333.2 |

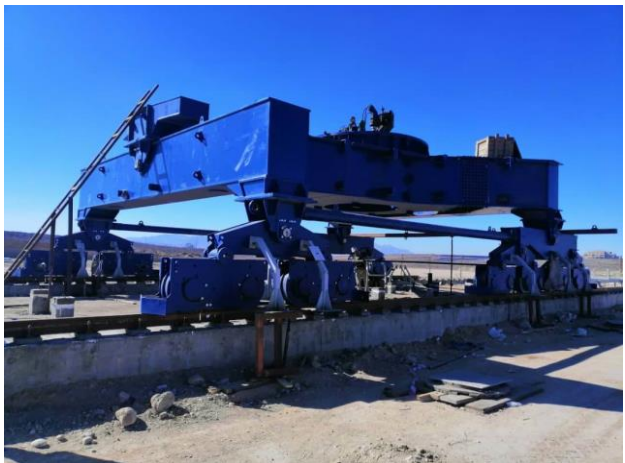

(a)

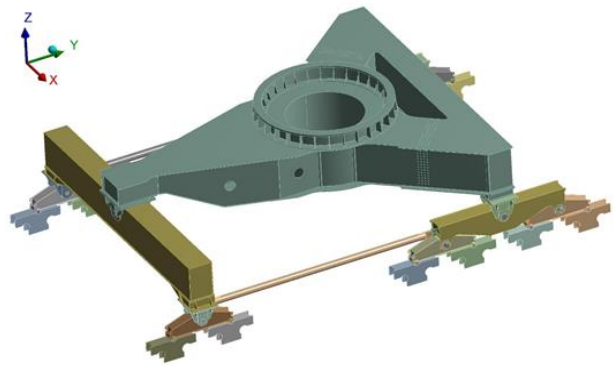

(b)

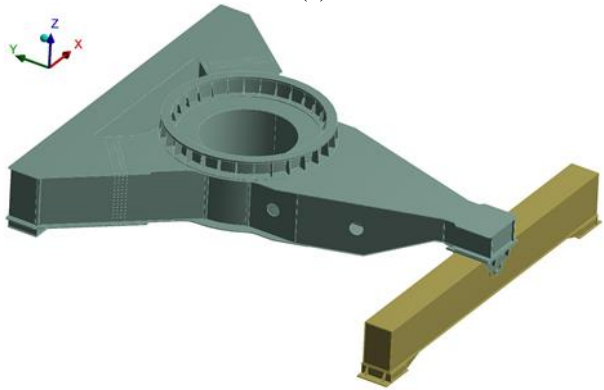

(c)

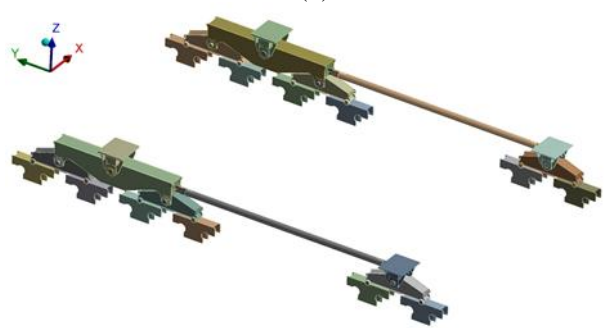

(d)

Figure 3. The undercarriage: a) built equipment, b) 3D model of the equipment, c) gantry and d) traveling system.

**Boundary Conditions**

In this analysis, the real movement of the wheels on the rail is simulated, see Fig. 4. Vertical displacement of all wheels is fixed. One wheel is fixed in all other directions. The other wheels are also constrained such that under loading they can move freely along the rails. Please note that, for simplification, the wheel assemblies are not modeled and the wheel boundary conditions are applied to the holes where they are attached to the traveling systems. Also, the joints which are shown in Fig. 5, are simulated using revolute-joints. A revolute joint is a one-degree-of-freedom kinematic pair used frequently in mechanisms and machines. The joints constrain the motion of two bodies to pure rotation along a common axis. Moreover, welded and bolted members are merged to each other.

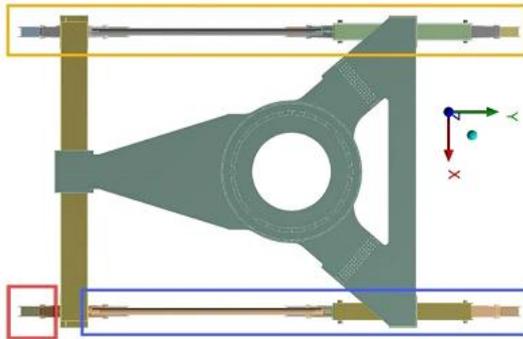

Figure 4. Wheels boundary conditions: Red: Dx,Dy,Dz fixed, Blue: Dy free Dx,Dz fixed, Yellow: Dy,Dx free Dz fixed.

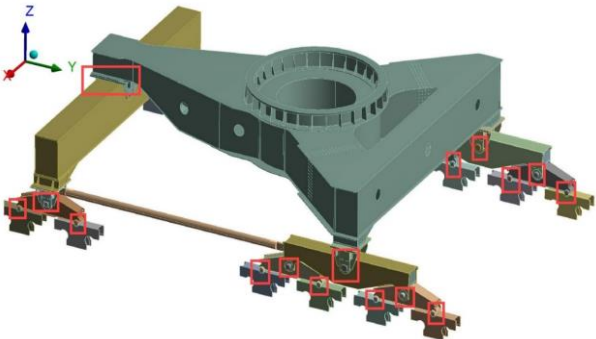

Figure 5. The revolute joints of the undercarriage (red).

**Loading**

There are several standards for different conditions of loading on stacker-reclaimer [4-10]. The forces and momentum on machine are calculated according to FEM2.131 [6] and ISO 5049-1 [5]. These forces include main loads, additional loads and special loads. These loadings are listed in Table 3.

The loads on the structure include the gravity force, the vertical force, as well as the moments applied to the location of the slew bearing (in reclaiming condition), see Fig. 6. The value and direction of slew bearing forces considering different working condition (idle, stacking and reclaiming condition) and angles of boom (slewing and luffing) are different. In this study, the most critical state, including reclaiming working, 45 degrees for slewing angle, and 9 degrees for luffing angle, have been considered. In this condition, the forces are as follows:

$$\begin{aligned} F_z &= -3{,}748 \text{ KN} \\ F_x &= F_y = 0 \text{ KN} \\ M_x &= 5{,}162 \text{ KN.m} \\ M_y &= -5{,}241 \text{ KN.m} \\ M_z &= 0 \text{ KN.m} \end{aligned} \quad (1)$$

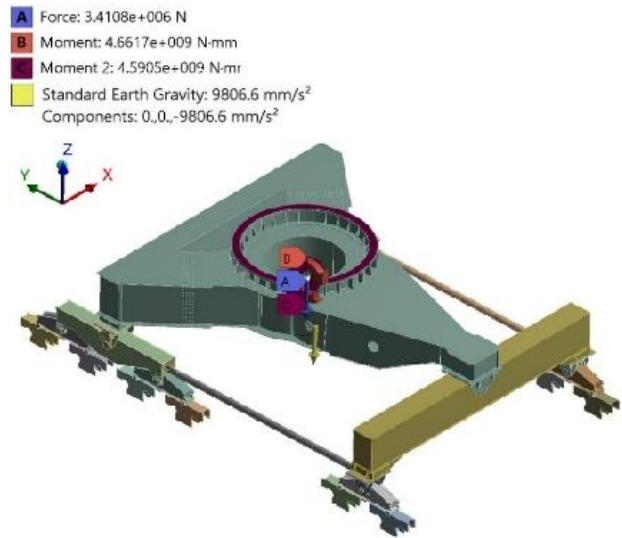

Figure 6. The gravity force and loading applied on top of the structure (slew bearing support).

**Finite element model**

In order to achieve the stress distribution, the geometry has been discretized by FEM package ANSYS. ANSYS-Workbench [10] is an interactive environment for simulation and solving all kinds of engineering and scientific problems based on partial differential equations (PDE). The whole geometry is discretized using approximately 120,000 tetrahedral 2D elements, see Fig. 7. Solid modeling is employed because of its simplicity and accuracy in comparison with shell modeling.

**Results and Discussion**

Von Mises effective stresses are evaluated and shown in Figs. 8. Deformation of the whole structure is shown in Fig. 9 magnified by a factor of 69 for better visualization of the deformed shapes. The maximum deflection is 14 mm which is in acceptable range. von Mises stress of the individual parts of the undercarriage structure is also shown in Fig. 10. The distribution of stresses over of all parts including gantry, connection structure and traveling system, is as follows. In the gantry, the highest stresses up to 190 MPa are observed on the top and bottom plates; at some points near to the slew bearing flange. In the traveling system, the wheels' supporting structures connected to the traveling system are in the low-risk stress state. In this system, there is four bridges to connect the upper and bottom structures, including secondary balancing beam and support structures. In these bridges, the maximum stress, which is located in the middle of structures is less than 100 MPa.

Table 3. loads according to FEM 2.131.

| Loading Categories | Loads | Location of load |
|---|---|---|
| Main Load | Dead load | Structure |
| | | Counterweight |
| | | Equipment |
| | Material Load | On conveyor |
| | | On bucket wheels |
| | | In the hopper |
| | Incrustation | On conveyor |
| | | On bucket wheels |
| | Normal digging force | Tangential |
| | | Lateral |
| Additional loads | Wind | In service wind |
| | Abnormal digging force | Tangential |
| | | Lateral |
| | Snow and ice load | Vertical |
| | Friction | Rail and wheels |
| | Skewing | Perpendicular to the rail by traveling |
| Special load | Wind | Out of service wind |
| | Seismic | Horizontal in different angles |
| | Buffer effect | On buffer |
| | Blocking of traveling device | On wheel |

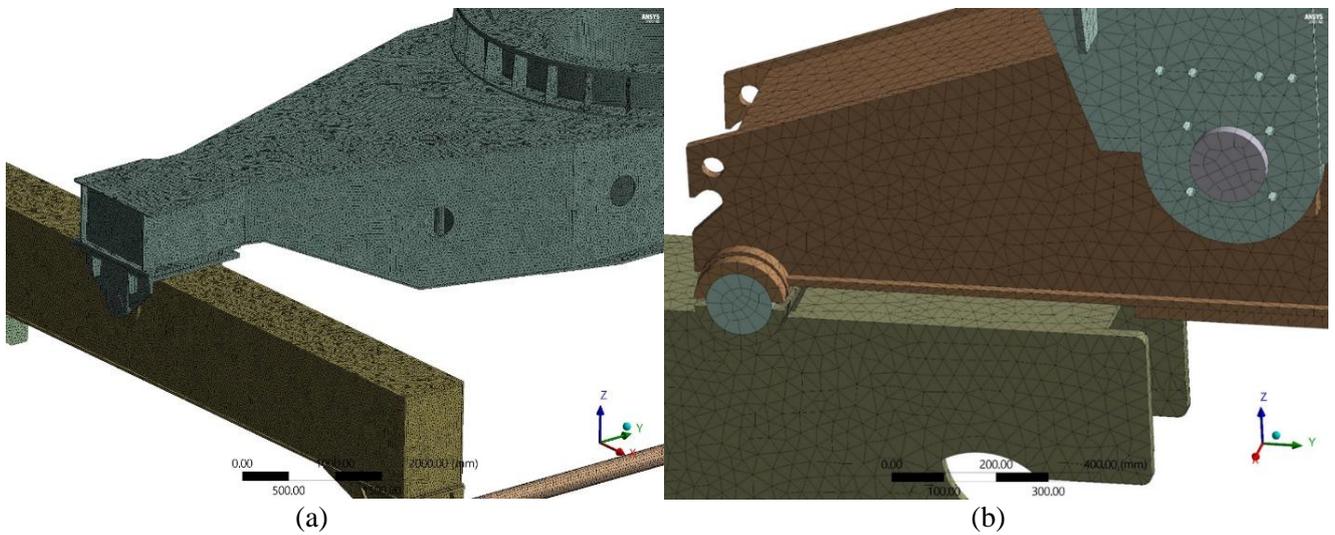

(a) (b)

Figure 7. Views on the finite element model.

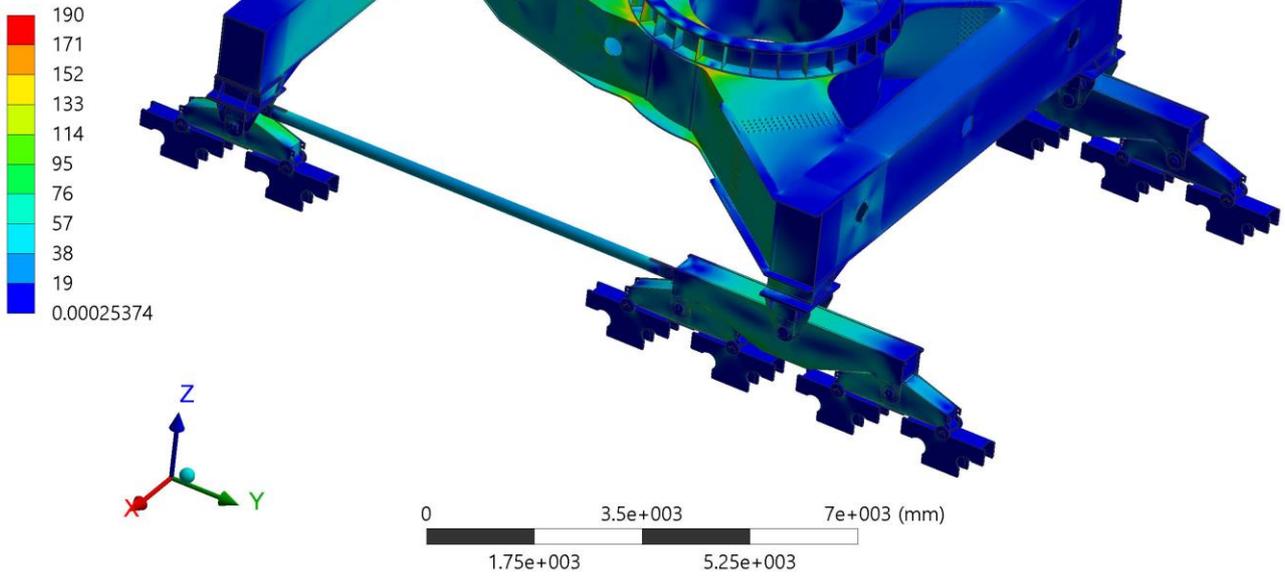

Figure 8. von Mises stress distribution on the undercarriage structure.

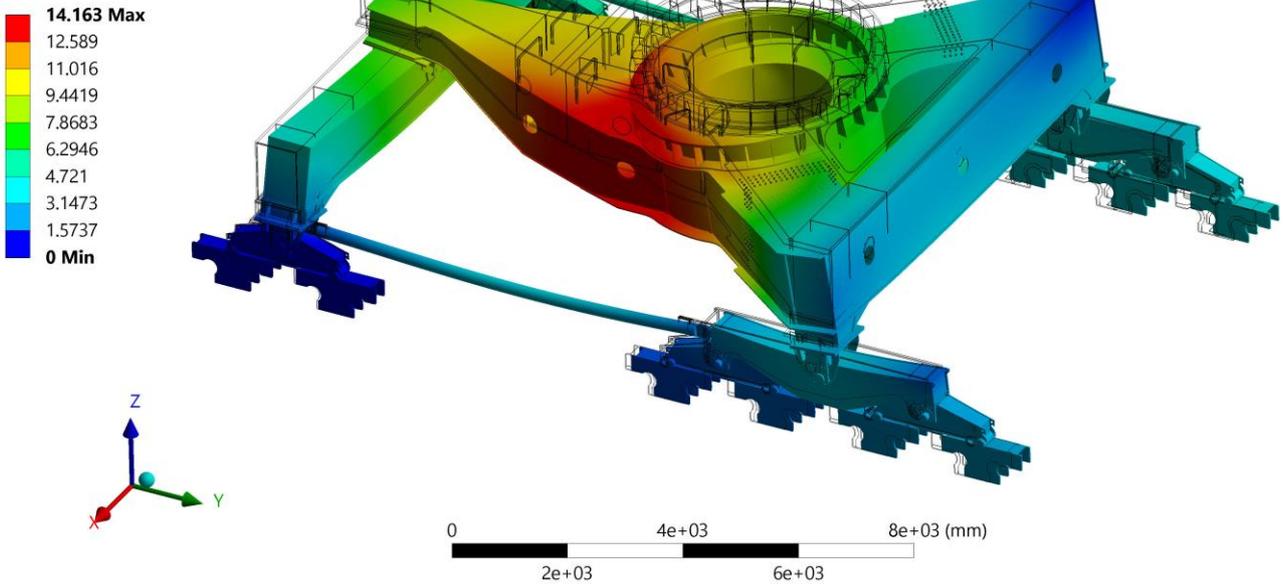

Figure 9. Deformation of the undercarriage structure (magnification factor: 69)

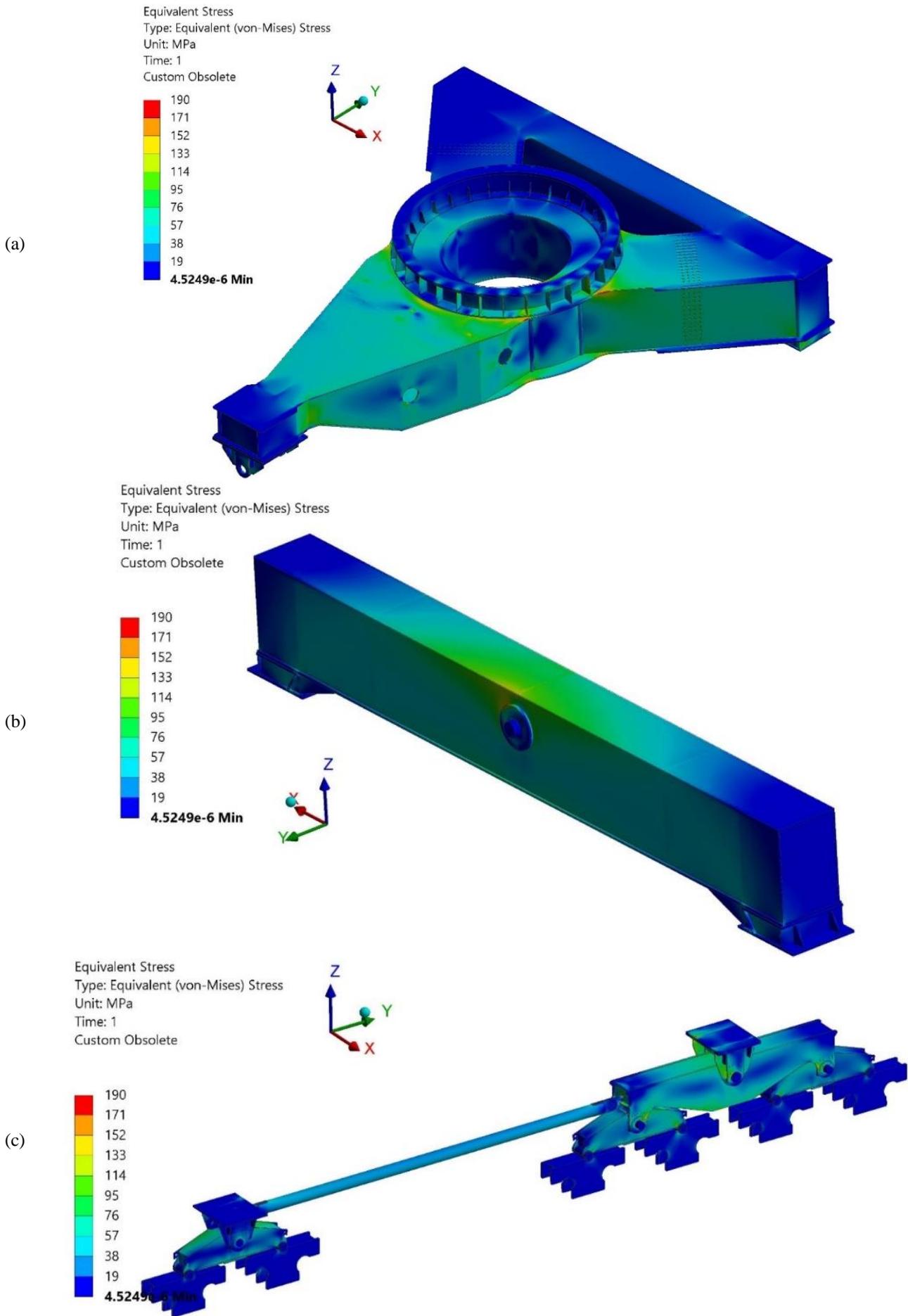

Figure 10. von Mises stress distribution: a) gantry, b) bridge on connection and c) traveling system.

## Conclusion

A combined stacker-reclaimer machine at Goharzamin Iron Ore Company which is being designed and constructed by a local supplier for the first time in Iran. It is highly recommended in the machine standards [5-9] particularly [9] to audit the design of the equipment in order to prevent machine failure. The main goal of the current study was to investigate static strength of the machine structure (undercarriage) under critical loading condition using Finite Element Analysis.

Based on the results, the following conclusions can be made:

- All welded parts are merged and the whole solid model (3D) of the structure was simulated with tetrahedral elements generated by free meshing. Compared with approaches such as shell and beam modeling, no further simplification and preprocessing is required and the obtained results are quite convincing. The structural joints and the wheels are also simply included in the model.
- The maximum deflection is 14 mm which is in acceptable range.
- In the gantry, the highest stresses up to 190 MPa are observed on the top and bottom plates; at some points near to the slew bearing flange.
- In the traveling system, the wheels' supporting structures connected to the traveling system are in the low-risk stress state.
- The maximum stress in bridges (for connection of upper structure to traveling system), which is located in the middle of structures is less than 100 MPa.

Design of the undercarriage structure was static analyzed and shown to be appropriate. In the following, the fatigue strength shall be investigated according the machine standards [5-9] and and the standards such as BS 7608 [15] for welded structures.
the static and fatigue strength of the super structure should also be audited with the same approach.

## Acknowledgement
The authors would like to express their special thanks for the provided supports by the Central Technical Office at Goharzamin Iron Ore Company (Sirjan).